\batchmode
\documentclass[twoside]{article}

\usepackage{amsbsy}
\usepackage{amssymb}
\usepackage{amsfonts}
\usepackage{amsmath}
\usepackage{amsxtra}
\usepackage{exscale}
\usepackage{enumerate}
\usepackage{hhline}
\usepackage{theorem}
\usepackage{array}
\makeatletter
\renewcommand{\section}{\@startsection
{section}%
{1}%
{0mm}%
{-1.5\baselineskip}%
{1.0\baselineskip}%
{\normalfont\large\bf}}%
\makeatother
 {\theoremstyle{change}
\theoremheaderfont{\normalfont\bfseries}
\theorembodyfont{\slshape}
\newtheorem{Theo}{Theorem:}[section]}
{\theoremstyle{change} \theorembodyfont{\slshape}
\newtheorem{Prop}[Theo]{Proposition:}}
{\theoremstyle{change} \theorembodyfont{\slshape}
}
{\theoremstyle{change} \theorembodyfont{\slshape}
\newtheorem{Cor}[Theo]{Corollary:}}
{\theoremstyle{change} \theorembodyfont{\upshape}
\newtheorem{Rem}[Theo]{Remark:}}
{\theoremstyle{change} \theorembodyfont{\upshape}
}
{\theoremstyle{change} \theorembodyfont{\upshape}
}
{\theoremstyle{change} \theorembodyfont{\upshape}
\newtheorem{Not}[Theo]{Notation:}}
{\theoremstyle{change} \theorembodyfont{\upshape}
\newtheorem{Nr}[Theo]{\kern -1ex}}
{\theoremstyle{change}\theorempreskipamount=0cm\theorempostskipamount=0cm
\theorembodyfont{\upshape\small}
}

\newcommand{\NN}{\mathbb{N}}

\DeclareMathOperator{\ord}{ord}

\DeclareMathOperator{\HN}{HN}
\begin{document}
\title{On the  value-semigroup  of a simple complete ideal in a two-dimensional regular local ring}
\author{S.\ Greco\\Politecnico di Torino\begin{footnote}{email: silvio.greco@mail.polito.it}\end{footnote} \and
K.\ Kiyek\\University of Paderborn\begin{footnote}{email: karlh@math.upb.de}\end{footnote}}
\date{}
\maketitle

\begin{abstract} {Let $R$ be a two-dimensional regular local ring with maximal ideal
$\mathfrak m$, and let $\wp$ be a simple complete $\mathfrak m$-primary ideal which is residually rational.
Let $R_0:=R\subsetneqq R_1\subsetneqq \cdots\subsetneqq R_r$ be the quadratic sequence associated to $\wp$, let
$\Gamma_\wp$ be the value-semigroup associated to $\wp$, and let $(e_j(\wp))_{0\leq j\leq r}$
be the multiplicity sequence of $\wp$. We associate to $\wp$ a sequence
$(\gamma_i(\wp))_{0\leq i\leq g}$ of natural integers which we call the formal characteristic sequence of $\wp$. We show
that the value-semigroup of $\wp$, the multiplicity sequence of $\wp$ and
the formal characteristic sequence of $\wp$
are equivalent data; furthermore, we give a new proof of the fact that $\Gamma_\wp$ is symmetric. In order to prove
these results, we use the Hamburger-Noether tableau of $\wp$ [\,cf.\ \cite{greki}\,]; we also give a formula for $c_\wp$, 
the conductor of $\Gamma_\wp$, in terms of the Hamburger-Noether tableau of $\wp$}\end{abstract}
\setcounter{section}{0}
\section{Introduction and Notations}
 \begin{Nr}
{\sc Introduction:} \label{seq0} Let $R$ be a two-dimensional
regular local ring with maximal ideal $\mathfrak m$, quotient field
$K$ and residue field $\kappa$, and let $\wp $ be a simple
complete $\mathfrak m$-primary ideal of $R$. We use the terminology
and  results developed in chapter VII of \cite{kivic}; for other
references cf.\ Appendix 3 of \cite{Zar} and the papers
\cite{Li4}, \cite{Li2} and \cite{Li3} of Lipman. The ideal $\wp$
determines a sequence $R_0:=R\subsetneqq
R_1\subsetneqq\cdots\subsetneqq R_r=:S$ where $R_1,\ldots,R_r$ are
two-dimensional regular local rings with quotient field $K$; for
each $j\in\{1,\ldots,r\}$ the ring $R_j$ is a quadratic transform
of $R_{j-1}$, and the residue field $\kappa_j$ of $R_j$ is a
finite extension of $\kappa$; we denote its degree by
$[\,R_j:R\,]$ [\,cf.\ \cite{kivic}, VII(5.1)\,]. We set
$f_\wp:=[\,R_r:R\,]$, and we call the ideal $\wp$ residually
rational if $f_\wp=1$. The order function $\ord_S$ of $S$ gives
rise to a discrete rank $1$ valuation of $K$ which shall be
denoted by $\nu_\wp=\nu$; let $V$ be the valuation ring of $\nu$.
Then, for $j\in\{1,\ldots,r\}$, $R_j$ is the only quadratic
transform of $R_{j-1}$ which is contained in $V$.

(1) The semigroup $\Gamma_\wp:=\{\nu(z)\mid z\in
R\setminus\{0\}\}$ is called the value-semigroup of $\wp$. Since
$K$ is the quotient field of $V$ and $R$, there exist non-zero
elements $z_1$, $z_2\in R$ with $\nu(z_1)-\nu(z_2)=1$, hence
$\Gamma_\wp$ is a numerical semigroup, i.e., the greatest common
divisor of its elements is $1$, and therefore $\Gamma_\wp$ is
finitely generated. The group $\Gamma_\wp$  was studied in papers
of Noh \cite{Noh1},
 Lipman \cite{Li3}, Spivakovsky \cite{Spiv} and in \cite{greki};
in \cite{Noh1} Noh showed---in case $\kappa$ is algebraically
closed---that $\Gamma_\wp$ is symmetric; Lipman \cite {Li3} showed
that this result also holds if $\wp$ is residually rational and he
gave a formula for the conductor $c_\wp$ of $\Gamma_\wp$ in terms
of the multiplicities of $\wp$ [\,cf.\ (2) for this notion\,]. In
the rest of this note, {\sl we always assume that $\wp$ is
residually rational.} We choose a system $\widetilde\kappa$ of
representatives of $\kappa$ in $R$ with the additional property
that the zero element of $\kappa$ is represented by $0$. Since all
the rings $R_1,\ldots,R_r$ have the same residue field $\kappa$,
we can take $\widetilde \kappa$ as a set of representatives of
$\kappa$ in each of the rings $R_1,\ldots,R_r$. In \cite{greki},
Prop.\ (8.6), we constructed---using the Hamburger-Noether tableau
of $\wp$---a system of generators $r_0,\ldots,r_h$ of $\Gamma_\wp$
such that $\Gamma_\wp$ is even strictly generated by these
elements [\,cf.\  \cite{greki}, (8.3),  for this notion\,].

(2) For $j\in\{0,\ldots,r\}$ let $\wp^{R_j}$ be the transform of
$\wp $ in $R_j$ [\,cf.\ \cite{kivic}, chapter VII, section
(2.3)\,]; note that $\wp^{R_r}$ is the maximal ideal of $R_r$.  We
set
$$
e_j(\wp)=\ord_{R_j}(\wp^{R_j})\quad\mbox{for
$j\in\{0,\ldots,r\}$.}
$$
The sequence $(e_j(\wp))_{0\leq j\leq r}$ is called the
multiplicity sequence of $\wp$ [\,these  integers are the non-zero
elements of  the point basis of $\wp$, cf.\ Lipman \cite{Li2},
\cite{Li3}\,].

 (3) Let $ n_0> n_1> \cdots> n_t = 1$ be the distinct multiplicities
  occurring in the
sequence $(e_j(\wp))_{0 \leq j\leq r}$, and assume that, for every
$ i \in \{1, \ldots, t\}$, $ n_i$ occurs exactly $h_i$ times. Then
the multiplicity sequence $(e_j(\wp))_{0 \leq j\leq r}$ of $\wp$
takes the form
$$
\underbrace{n_0, \ldots, n_0}_{h_0},\ \underbrace{n_1, \ldots,
n_1}_{h_1}, \ldots, \underbrace{n_{t-1}, \ldots,
n_{t-1}}_{h_{t-1}}, \,\,\underbrace{1,1, \ldots, 1}_{h_t}.
$$
We call
$$
\underbrace{n_0, \ldots, n_0}_{h_0},\ \underbrace{n_1, \ldots,
n_1}_{h_1}, \ldots, \underbrace{n_{t-1}, \ldots,
n_{t-1}}_{h_{t-1}}
$$
the essential part of the multiplicity sequence; note that
$$
h_t=\ell_R(R/\wp)-\sum_{i=0}^{t-1}\frac{h_in_i(n_i+1)}{2}
$$
by the length formula of Hoskin-Deligne [\,cf.\ \cite{kivic},
chapter VII, section (6.9)\,].

 Furthermore, let
$$
 \{s_1, \ldots, s_g\} = \{j \in \{1,
\ldots, t\} \mid  n_j \mbox{ divides } n_{j-1}\},
$$
and let these elements be labeled in such a way that $s_1 < \cdots
< s_g$. Note that $ g=0$ iff $e_0(\wp) = 1$.  The non-negative
integer $g$ will be called the genus of $\wp$. It is often
convenient to set $s_0 = 0$. We set
$$
f_i:=\frac{n_{s_i-1}}{n_{s_i}}\quad\mbox{for
$i\in\{1,\ldots,g\}$.}
$$

(4) We set $\gamma_0(\wp):=n_0$, and we define, by
recursion\begin{footnote}{the elements $h_0, h_1\ldots$ here are
the elements which are denoted by $e_0,e_1,\ldots$ in \cite{Ki8},
(2.2) and (2.3).}\end{footnote}
\begin{align*}\tag{$*$}
\gamma_1(\wp)&:=h_0n_0+n_1,\\\tag{$**$} \gamma_{i+1}(\wp)
-\gamma_i(\wp) &:= (h_{s_i} - f_i) n_{s_i} + n_{s_i+1} \quad
\mbox{for every } i \in \{1, \ldots, g-1\}.
\end{align*}
We call the sequence $(\gamma_i(\wp))_{0\leq i\leq g}$ the formal
characteristic sequence of $\wp$.

(5) We show in this note,  using only the Hamburger-Noether
tableau of $\wp$, that $\Gamma_\wp$ is symmetric, give several
formulae which express $c_\wp$, the conductor of $\Gamma_\wp$, in
terms of data of the Hamburger-Noether tableau of $\wp$, give a
new proof for the formula of Lipman \ref{lipman} expressing
$c_\wp$ in terms of the multiplicities of $\wp$, and show that
$\Gamma_\wp$, the sequence of multiplicities of $\wp$, and the
formal characteristic sequence of $\wp$ are equivalent data.
\end{Nr}

\begin{Nr}\label{seq1}{\sc Multiplicity sequence and quadratic transforms:}  Let
$i\in\{0,\ldots,t-1\}$,  set $j:=h_0n_0+\cdots+h_in_i$ and
$T:=R_j$. Let $T=T_0,T_1,\ldots$ be the   sequence of quadratic
transforms of $T$ along $V$; we have $n_i=e(\wp^T)$,
$n_{i+1}=e(\wp^{T_1})$. There exists a regular system of
parameters $\{x,y\}$ for $T$ such that $\nu(x)=n_i$,
$\nu(y)>\nu(x)$, and $\nu(y/x)=n_{i+1}$.

(1) Firstly, we consider the case that $n_{i+1}$ does not divide
$n_i$. Then we have  $n_{i+1}=e(\wp^{T_k})$ for
$k\in\{1,\ldots,h_{i+1}\}$, and, setting $x_1:=y/x$, $y_1:=x$,
$\{x_1,y_1\}$ is a regular system of parameters for $T_1$, and
since $\nu(y_1/x_1^{k-1}) = n_i-(k-1)n_{i+1}$ for
$k\in\{1,\ldots,h_{i+1}\}$, it is clear that  $\{x_k:=x_1,
y_k:=y_1/x_1^{k-1}\}$ is a regular system of parameters for $T_k$
for $k\in\{1,\ldots,h_{i+1}\}$. We have
$\nu(y_{h_{i+1}})=n_i-(h_{i+1}-1)n_{i+1}\geq n_i$,
$n_i-h_{i+1}n_{i+1}<n_i$, hence $n_i-h_{i+1}n_{i+1}=n_{i+2}$, and
this is the division identity for the natural integers $n_i$,
$n_{i+1}$.

(2) Now, we treat the case that $n_{i+1}$ divides $n_i$, and set
$f:=n_i/n_{i+1}$. We consider the first $f$ quadratic transforms
$T_1,\ldots,T_f$ of $T$ along $V$. Setting $x_k:=x_1:=y/x$,
$y_k:=y_1/x_1^{k-1}$ for $k\in\{1,\ldots,f\}$, we see that
$\{x_k,y_k\}$ is a regular system of parameters for $T_k$, and
that $T_{f+1}$, the quadratic transform of $T_f$ along $V$, has a
system of parameters $\{x_1, (y_1-a x_1^f)/x_1^f\}$ where $a\in R$
is the unique element in $\widetilde\kappa\setminus\{0\}$ such
that $\nu((y_1-ax_1^f)/x_1^f)>0$ [\,for the existence of such an
element $a$ remember that $\wp$ is residually rational, that
$T_0\neq R_r$, and use the argument given in \cite{greki},
(7.5)(3)\,], and $e(\wp^{T_f})\geq e(\wp^{T_{f+1}})$. This
implies, in particular, that $f\leq h_{i+1}$.
\end{Nr}

\begin{Nr}
{\sc Multiplicity sequence and formal characteristic sequence:} We
show that the formal characteristic sequence
$(\gamma_i(\wp))_{0\leq i\leq g}$ of $\wp$ determines the
essential part of the multiplicity sequence $(e_j(\wp))_{0\leq
j\leq r}$ of $\wp$, hence that the multiplicity sequence and the
formal characteristic sequence of $\wp$ are equivalent data.

(1) We have by \ref{seq1}
\begin{align*}
n_{j-1} &= h_j n_j + n_{j+1}&&\quad\mbox{for }j \in \{1, \ldots,
r-1\} \setminus \{s_1, \ldots,
s_g\}, \\
n_{s_i-1}& =  f_i n_{s_i}\mbox{ where } f_i \in \NN \mbox{ and }
f_i \leq h_{s_i} &&\quad\mbox{for } i \in \{1, \ldots, g\}.
\end{align*}
(2) Let $i \in \{0, \ldots, g-1\}$. The continued fraction
expansion of $n_{s_i}/n_{s_i+1}$ takes the form
$$
\frac{n_{s_i}}{n_{s_i+1}} = [\,h_{s_i+1}, \ldots, h_{s_{i+1}-1},
f_{i+1}\,],
$$
 and we have
$$\kern-4.63pt\hbox{$
\gcd (n_{s_i}, n_{s_i+1}) = \gcd (n_{s_i}, n_{s_i+1}, n_{s_i+2}) =
\cdots = \gcd (n_{s_i}, n_{s_i+1}, \ldots, n_{s_{i+1}-1}) =
n_{s_{i+1}}.$}
$$
From this it follows easily that
$$
\gcd (n_0, n_1, \ldots, n_{s_i-1}) = n_{s_i} \quad \mbox{for every
} i \in \{1, \ldots, g\}.
$$
(3) Let $\gamma_i:=\gamma_i(\wp)$ for $i\in\{0,\ldots,g\}$. We
show that
\begin{align*}
\gamma_0&=n_0,\\
\gcd (\gamma_0, \ldots, \gamma_i) &= \gcd (n_{s_{i-1}},
n_{s_{i-1}+1}) = n_{s_i} \quad \mbox{for every } i \in \{1,
\ldots, g\}.
\end{align*}
This holds if $g=0$. Now we assume that $g\geq1$.  The assertion
holds if $ i = 1$ since, by definition,  $ \gamma_0 = n_0$, $
\gamma_1 = h_0n_0+n_1$, hence
$\gcd(\gamma_0,\gamma_1)=\gcd(n_0,n_1)=n_{s_1}$ [\,cf.\ (2)\,].
Assume that $ i \in \{1, \ldots, g-1\}$, and that the assertion
holds for $i$, hence, in particular, that $n_{s_i}$ divides
$\gamma_i$. Then $\gcd (\gamma_0, \ldots, \gamma_{i+1}) = \gcd
(n_{s_i}, \gamma_{i+1})$. From (2) and \ref{seq0}(4)$(**)$ we see
that  $ n_{s_{i+1}} = \gcd (n_{s_i}, \gamma_{i+1}) = \gcd
(n_{s_i}, n_{s_i+1})$.

(4)  Now $(*)$ in \ref{seq0}(4) is the division identity for the
integers $\gamma_1$ and $\gamma_0=n_0$. By (3) we have $ n_{s_i} =
\gcd (\gamma_0, \ldots, \gamma_i)$ for every $ i \in \{0, \ldots,
g\}$. Let $ i \in \{1, \ldots, g-1\}$; since $ n_{s_i}
> n_{s_i+1}$, $(**)$ in \ref{seq0}(4) is the
division identity for the integers $ \gamma_{i+1} - \gamma_i$ and
$n_{s_i}$. Knowing $ n_{s_i}$ and $ n_{s_i+1}$ for every $ i \in
\{0, \ldots, g-1\}$, we use (2) to calculate $ n_0, \ldots, n_r$,
$ h_0, \ldots, h_{r-1}$, i.e., we have calculated the essential
part of the multiplicity sequence of $\wp$ from its formal
characteristic sequence.

\end{Nr}
\section{Quadratic Transform and $\HN$-Tableau}
\begin{Nr} {\sc Quadratic transform:}
 We consider the quadratic transform $R_1$ of $R$ and let
$n'_0>n'_1>\cdots>n'_{t'}$ etc., be the numbers defined by the
multiplicity  sequence $(e(\wp^{R_j}))_{1\leq j\leq r}$ of the
simple ideal $\wp^{R_1}$ of $R_1$. It is easy to check that the
connection  between the sequence  $n_0,\ldots,n_t$ and the
sequence  $n'_0,\ldots,n'_{t'}$  is given by the  formulae  in
\cite{Ki8}, (2.4)-(2.6).
\end{Nr}

\begin{Nr}
\label{nohgroup1} {\sc Hamburger-Noether tableau and quadratic
transform:}  Let $\{x,y\}$ be a regular system of parameters for
$R$. We assume that $\wp\neq\mathfrak m$;  remember that $\nu$
dominates $R_1$. If $\nu(x)<\nu(y)$, then we take $\{x_1:=x,
y_1:=y/x\}$ as a regular system of parameters for $R_1$, if
$\nu(x)>\nu (y)$, then we take $\{x_1:=y,y_1:=x/y\}$ as a regular
system of parameters for $R_1$; finally,   if $\nu(x)=\nu(y)$,
then let $a_1\in\widetilde \kappa\setminus \{0\}$ be the unique
element with $\nu(y-a_1x)>\nu(x)$ [\,for the existence of $a_1$,
note that $r\geq1$, and   argue as in \ref{seq1}(2)\,], and in
this case we take $\{x_1:=x, y_1:=(y-a_1x)/x\}$ as a regular
system of parameters for $R_1$.

Let
$$
\HN(\wp;x,y) = \begin{pmatrix}  p_i\\c_i\\a_i \end{pmatrix} _{1
\le i \leq l}
$$
be the Hamburger-Noether tableau of $\wp$ with respect to
$\{x,y\}$ [\,cf.\ \cite{greki}, (7.6), for the definition of the
HN-tableau $\HN(\wp;x,y)$\,]. We set $\wp':=\wp^{R_1}$, and we
assume that $\wp'$ is not the maximal ideal of $R_1$. Then the
Hamburger-Noether tableau of $\wp'$  with respect to the regular
system of parameters $\{x_1,y_1\}$ for $R_1$ constructed in the
preceding paragraph is
$$
\HN(\wp';x_1,y_1)= \begin{cases}
\begin{pmatrix}
p_1 - c_1 & p_2 & \cdots \\
c_1 & c_2 & \cdots \\
a_1 & a_2 & \cdots \end{pmatrix} &\mbox{if } p_1 > c_1,
\\[4ex]
\begin{pmatrix}
c_1 - p_1 & p_2 & \cdots \\
p_1 & c_2 & \cdots \\
a_1 &a_2 & \cdots
\end{pmatrix} &\mbox{if } p_1 < c_1, \\[4ex]
\begin{pmatrix}
p_2 & p_3 & \cdots \\
c_2 & c_3 & \cdots \\
a_2 & a_3 & \cdots \end{pmatrix}& \mbox{if } p_1 = c_1.
\end{cases}
$$

\end{Nr}

\begin{Rem}
\label{cahrtrans} Let $(\theta_1,q_1,\ldots,q_h)$ (resp.\
$(\theta'_1,q'_1,\ldots,q'_{h'})$) be the characteristic
sequence\begin{footnote}{Cf.\  \cite{greki}, (8.4), for this
notion.}\end{footnote} of $\wp$ with respect to $\{x,y\}$ (resp.\
of $\wp'$ with respect to $\{x_1,y_1\}$). It is easy to check that
the connection between these  two sequences is given by the
formulae in \cite{Ki8}, (3.1)-(3.6) .

\end{Rem}

\begin{Prop}\label{III(8.19.1)}  Let $\{x,y\}$ be a
regular system of parameters for $R$,
 let $(\theta_1;q_1,\ldots,q_h)$ be the
characteristic sequence of\/ $\wp$ with respect to $\{x,y\}$, and
let $(\gamma_0,\ldots,\gamma_g)$ be the formal characteristic
sequence of $\wp$. Then we have \vspace{-2.5ex}
\begin{alignat*}{5}
 \intertext{\textup{(1)} if $ \theta_1 \leq q_1$ and  $ \theta_1 \nmid q_1$,
then}
 h &= g, & \gamma_0 &= \theta_1,&  \gamma_1&=q_1,&
\gamma_j-\gamma_{j-1}&=q_j &&\quad\mbox{for every } j\in \{ 2,
\ldots, g\},\\
\intertext{\textup{(2)} if  $ \theta_1 \mid q_1$ or $ q_1 \mid
\theta_1$, then} h &= g+1, &\quad\gamma_0&=\theta_1,&\quad
\gamma_1&= q_1+q_2, & \quad\gamma_i-\gamma_{i-1}&= q_{i+1}&&
\quad\mbox{for every }j \in\{ 2,
\ldots, g\},\\
\intertext{\textup{(3)} if $ q_1 < \theta_1$ and  $ q_1 \nmid
\theta_1$, then} h &= g, & \gamma_0&= q_1, & \gamma_1 &=
\theta_1,& \gamma_j-\gamma_{j-1}&= q_j &&\quad\mbox{for every } j
\in\{2, \ldots, g\}.
\end{alignat*}

\end{Prop}

Proof: By induction, we may assume that the corresponding formulae
hold for the $\HN$-tableau $\HN(\wp';x_1,y_1)$, defined as in
\ref{nohgroup1}; it is a routine matter, using \ref{seq1} and
\ref{cahrtrans}, to verify these relations.

\begin{Cor}
\label{endresult} Let $\wp$ be a simple complete  $\mathfrak
m$-primary
 ideal of $R$ which is residually
rational. The following data are equivalent:

\textup{(1)} The semigroup $\Gamma_\wp$ of $\wp$.

\textup{(2)} The multiplicity sequence of $\wp$.

\textup{(3)} The formal characteristic sequence of $\wp$.
\end{Cor}

Proof: This is an easy consequence of the last results.
\begin{Rem} (1) For quite a different proof of the fact
 that the data in (1) and (2) of \ref{endresult} are equivalent
 if $\kappa$ is algebraically closed, cf.\ the proof of
Th.~3 in Noh's paper \cite{Noh1}.

(2) Let $\kappa$ be an algebraically closed field of
characteristic $0$, let $\kappa[\![\,t]\!]$  be the ring of formal
power series over $\kappa$, let $(\gamma_i)_{0\leq i\leq g}$ be
the formal characteristic sequence of $\wp$, and set $
 x:=t^{\gamma_0}$,\ $y:=t^{\gamma_1}+\cdots+t^{\gamma_g}$;
 let $C$ be the plane algebroid curve defined by this
 parametrization. Then $(\gamma_i)_{0\leq i\leq g}$ is the
 characteristic sequence of this curve with respect to the pair $(x,y)$. If one determines the HN-tableau for the pair $(x,y)$
 [\,cf.\ \cite{Russell}, (2.2)\,], and calculates its semigroup
 sequence [\,cf.\ \cite{Russell}, 2.8.5\,], then it is easy to check
 that the semigroup which is strictly generated by this semigroup
 sequence,
is the semigroup of $C$. This explains the name ``formal
 characteristic sequence''.
\end{Rem}
\begin{Prop}\label{symm} Let $\{x,y\}$ be a regular system of
parameters  for  $R$, let $\HN(\wp;x,y)$ be the $\HN$-tableau of
$\wp$ with respect to $\{x,y\}$, and let $(\theta_i)_{1\leq i\leq
h+1}$ (resp.\ $(r_i)_{0\leq i\leq h}$)  be the divisor sequence
 (resp.\ the semigroup sequence\begin{footnote}{Cf.\ \cite{greki},
 (8.4)(2) and (8.4)(3) for this notions.}\end{footnote})   of  $\wp$ with respect
to $\{x,y\}$. The semigroup $\Gamma_\wp$ is strictly generated by
the semigroup sequence $(r_i)_{0\leq i\leq h}$, it  is symmetric,
and its conductor $c_\wp$ is
$$
c_\wp= -r_0+1+\sum^{h}_{j=1} r_j(\theta_j/\theta_{j+1}-1).
$$
\end{Prop}

Proof: The assertion that $\Gamma_\wp$ is strictly generated by
the sequence $(r_i)_{0\leq i\leq h}$, was proved in \cite{greki},
Prop.~8.6; all the other statements  follow from the argument in
Russel's paper \cite{Russell}, proof of Th.~6.1, (2)-(6), where he
essentially uses only the fact that $\Gamma_\wp$ is strictly
generated by the sequence $(r_i)_{0\leq i\leq h}$.

\begin{Rem}
We take the opportunity to correct an omission in \cite{greki},
(8.17) and  (8.18).  We must add also the condition
$\theta_1>\cdots>\theta_{g+1}=1$ (this condition is used in order
to show that the tableau given in the proof of  \cite{greki},
(8.19) is, in fact, a HN-tableau). Here we should mention also the
paper of Angerm\"uller \cite{anger}, and, in particular, Satz 1.
\end{Rem}
\begin{Not} Let $\{x,y\}$ be a regular system of parameters for $R$, and
let
$$
\HN(\wp;x,y) = \begin{pmatrix}  p_i\\c_i\\a_i \end{pmatrix} _{1
\le i \leq l}
$$
be the $\HN$-tableau of $\wp$ with respect to $\{x,y\}$. We set
$$
\epsilon : = \min (\{i\in \NN_0 \mid c_j = c_{i+1} \mbox{ for
every $ j \geq i+1$\}}).
$$
\end{Not}
\begin{Cor}
\label{simpleidealAbh}  We keep the notations of
\textup{\ref{symm}}, and let $(\theta_1;q_1,\ldots,q_h)$ be the
characteristic sequence of $\wp$ with respect to $\{x,y\}$. For
the conductor $ c_\wp$ of the symmetric semigroup\/ $\Gamma_\wp$
we have
\begin{align*}
c_\wp & = {} -r_0+1+\sum^{h}_{j=1} r_j(\theta_j/\theta_{j+1}-1) \\
& = (q_1-1) (\theta_1-1) + \sum^{h}_{j=2} q_j (\theta_j - 1) \\
& =  (p_1-1) (c_1-1) + \sum^{\epsilon}_{i=2} p_i (c_i-1).
\end{align*}
\end{Cor}

Proof: We have shown in \ref{symm} that $\Gamma_\wp$ is symmetric,
and that the first expression in the above displayed formulae is
the conductor of $\wp$. From the definition of $\epsilon$ it
follows immediately that the second and third expression are
equal, and it is a trivial matter of arithmetic to show that the
first and the second expression are equal.

\begin{Cor}\textup{[\,Lipman\,]}\label{lipman}
We have
$$
c_\wp=\sum_{i=0}^re_i(\wp)(e_i(\wp)-1).
$$
\end{Cor}
Proof: Let $\{x_1,y_1\}$ be the regular system of parameters for
$R_1$ as constructed in \ref{nohgroup1}. We consider the second
expression for the conductor given in \ref{simpleidealAbh}, and
compare it  with the corresponding expression in
$\HN(\wp^{R_1};x_1,y_1)$; the difference equals
$e_0(\wp)(e_0(\wp)-1)$.


{\small
\begin{thebibliography}{99}
\bibitem{anger} Angerm\"uller, G., Die Wertehalbgruppe einer
ebenen irreduziblen algebroiden Kurve. Math.\ Z.\ {\bf 153}
(1977), 267--282.
\bibitem{cam} Campillo, A., Algebroid curves in positive
characteristic. LNM {\bf 813}. Springer, Berlin - Heidelberg - New
York 1980.
\bibitem{greki}  Greco, S.\ and Kiyek, K.,
General Elements of Complete Ideals and Valuations
   Centered at a Two-dimensional Regular Local Ring. In:
{A}lgebra, {A}rithmetic and {G}eometry with Applications,
381--455. Springer, Berlin - Heidelberg - New York 2004.

\bibitem{Ki8}
  Kiyek, K.,
       Multiplicity sequence and value semigroup.
Manuscripta Math.\ {\bf37}
 (1982), 211--216.
\bibitem{kivic} Kiyek, K.\ and Vicente J.\ L., Resolution of Curve
and Surface Singularities in Characteristic Zero. Kluwer,
Dordrecht 2004.
\bibitem{Li4}
 Lipman, J.,
    On complete ideals in regular local rings. In:
Algebraic geometry and commutative algebra, Vol.\ I,
   203--231. Kinokuniya, Tokyo 1988.
\bibitem{Li2}
   Lipman, J.,
  Proximity inequalities for complete ideals in two-dimensional
              regular local rings. In:
Commutative algebra: Syzygies, multiplicities, and birational
              algebra (South Hadley, MA, 1992),
   Contemp. Math.
  {\bf 159} (1992),
    293--306.


\bibitem{Li3}
 Lipman, J.,
         Adjoints and polars of
  simple complete ideals in two-dimensional regular local rings.
  Bull.\ Soc.\ Math.\ Belgique {\bf 45}
(1993),
  223--244.


\bibitem{Noh1}
       Noh, S.,
The Value Semigroup of Prime Divisors of
  the Second Kind in 2-dimensional Regular Local Rings.
 Trans.\ Amer.\ Math.\ Soc.\ {\bf336} (1993), 607--619.

\bibitem{Russell}
Russell, P.,
  Hamburger-{N}oether expansions and approximate
  roots of polynomials.
Manuscripta Math. {\bf31} (1980), 25--95.

\bibitem{Spiv} Spivakovsky, M., Valuations in Function Fields of
Surfaces. Amer.\ J.\ Math.\ {\bf112} (1990), 107--156.
\bibitem{Zar} Zariski, O.\ and Samuel, P., Commutative Algebra,
vol.\ II. Van Nostrand, Princeton 1960.
\end{thebibliography}
}
\end{document}